\documentclass[12pt,a4paper]{article}

\usepackage{amsmath,amsfonts}
\usepackage{mathrsfs} 
\usepackage{parskip}
\usepackage[T1]{fontenc}
\usepackage[utf8]{inputenc}
\usepackage{authblk}
\usepackage{enumerate}
\usepackage{graphicx}

\newtheorem{theorem}{Theorem}
\newtheorem{lemma}{Lemma}

\def\Pb{\mathbf{P}}

\def\Ex{\mathbf{E}}

\def\UU{\mathbb{U}}

\def\1{\mbox{1\hspace{-.25em}I}}
\begin{document}
\title{On Frequency Estimation   for Partially Observed    Processes with
  Small Noise in Observations. }
\author[1]{O.V. Chernoyarov}
\author[2]{Yu.A. Kutoyants}
\affil[1,2]{\small National Research University ``MPEI'', Moscow, Russia, }
\affil[2]{\small Le Mans University,  Le Mans,  France,}
\affil[1,2]{\small Tomsk State University, Tomsk, Russia}

\date{}

\maketitle
\begin{abstract}
We consider the problem of frequency estimation of the periodic signal
multiplied by a stationary Gaussian process (Ornstein-Uhlenbeck) and observed 
in the presence of the  white Gaussian noise. We show the consistency and
  asymptotic normality of the maximum likelihood  estimator  in the
asymptotics of small noise in observations. The model of observations is a linear
nonhomogeneous partially observed system and the construction and study of the
estimator  is essentialy based on the asymptotics of the equations of Kalman-Bucy filtration. 
\end{abstract}


\section{Introduction}
 The problem
of frequency estimation is of special interest in telecommunication theory.
The shift of the frequency allows to estimate the speed of the object
(Doppler  effect). The estimation of the frequency of the Gaussian signals
observed in the WGN like the studied in this work is of interest in
statistical radio physics.

This is the second work devoted to the frequency estimation by the
observations of partially observed linear nonhomogeneous system.
 In the first work we studied the properties of the
maximum likelihood estimators (MLE) and Bayes estimators (BE) in the situation,
where the noises in the state and observation equations tend to zero
\cite{ChK1}. It was shown that these estimators are consistent and
asymptotically normal. Here we study a slightly different model of
observations, where the noise in the state equation does not tend to zero and
we have asymptotics of {\it small noise}
  in the observations equation. The detailed study of dynamical
systems with small perturbations (noise) can be found in \cite{FV}. The
statistical problems for such models are presented in \cite{Kut94}.

Let us remind the models of
observations.  Consider the following model of observations
\begin{align}
\label{0}
 {\rm d}X\left(t\right)=A\cos\left(2\pi \vartheta t\right)Y_t{\rm
   d}t+\sigma{\rm d}W_t ,\quad X_0=0, \qquad 0\leq t\leq T 
\end{align}
where $Y_t,0\leq t\leq T$ is some stationary Gaussian process with
the spectral density function $f\left(\lambda \right)$ and
$W_t,0\leq t\leq T$ is the Wiener process. We have to
estimate the frequency $\vartheta $ by the continuous time observations
$X^T=\left(X_t,0\leq t\leq T\right)$.  
We suppose that the process $Y_t$ satisfies a linear
equation. The most simple examples are given below.

{\bf Example 1.} Ornstein - Uhlenbeck (O-U) process
\begin{align*}
{\rm d} Y_t=-a Y_t{\rm d}t+b {\rm d} V_t,\quad Y_0=y_0,\qquad 0\leq t\leq T
\end{align*}
has the correlation function $R\left(\tau \right)$ and the 
spectral density $f\left(\lambda \right)$  are the following
\begin{align*}
R\left(\tau \right)=\frac{b^2}{2a}e^{-a\left|\tau \right|},\qquad \quad
f\left(\lambda \right)=\frac{b^2}{a^2+4\pi ^2\lambda ^2}. 
\end{align*}
 We suppose that $a>0,b>0$ and $V_t$ is a Wiener process independent of $W_t$.

{\bf Example 2.} Suppose that the process $Y_t$ satisfies the equation
\begin{align*}
{\rm d}\dot Y_t=-a_1\dot Y_t{\rm d}t- a_0Y_t{\rm d}t+b  {\rm d}V_t,\qquad 0\leq t\leq T,
\end{align*}
where $a_1^2<4a_0$. Then it has the correlation function
\begin{align*}
R\left(\tau \right)=\frac{b^2}{2a_0a_1}e^{-\alpha \left|\tau \right|}\left(
\cos\left(\beta \tau \right)+\frac{\alpha }{\beta }\sin\left(\beta \tau
\right)  \right) , 
\end{align*}
where $\alpha =a_1/2$ and $\alpha ^2+\beta ^2=a_0$. 
 The spectral density is
\begin{align*}
f\left(\lambda \right)=\frac{b^2}{ \left[a_0-4\pi ^2\lambda ^2\right]^2+4\pi
  ^2a_1^2\lambda ^2} .
\end{align*}

 In this  work like \cite{ChK1} we suppose that $Y_t$ is O-U
process and we use the formalism of stochastic calculus. Note that all
results presented in this work for O-U process can be directly extended on the
model of Example 2.

We consider a slightly more general model, where the signal $A\cos\left(2\pi
\vartheta t\right)$ we replace by a known smooth periodic function
$f\left(\vartheta t\right)$. For simplicity we suppose that the period is
equal 1.  Of course, the signal $f\left(\vartheta t\right)$ has period $\tau
=\frac{1}{\vartheta }$.  Therefore we have a two-dimensional stochastic
process $\left(X_t,Y_t,0\leq t\leq T\right)$ satisfying the differential
equations
\begin{align}
\label{a}
{\rm d}X_t&=f\left(\vartheta t\right)Y_t{\rm d}t+\sigma {\rm
  d}W_t,\qquad X_0=0, \\ 
{\rm d}Y_t&=-aY_t{\rm d}t+b {\rm
  d}V_t,\qquad\qquad Y_0=y_0.
\label{b}
\end{align}
Here $W_t,0\leq t\leq T$ and $V_t,0\leq t\leq T$ are two independent Wiener
processes. The parameters $A, \sigma , a, b$ are supposed to be known and
positive, the parameter $\vartheta \in \Theta =\left(\alpha ,\beta \right)$
(frequency) is unknown and has to be estimated by the observations
$X^T=\left(X_t,0\leq t\leq T\right) $.

We study the maximum likelihood estimator (MLE) $\hat\vartheta _T$.  Let us
remind the definition of it. Introduce the conditional
mathematical expectation
\begin{align*}
m\left(\vartheta,t \right)=\Ex_\vartheta \left(Y_t\,|\,X_s,0\leq s\leq
t\right),\qquad 0\leq t\leq T .
\end{align*}
The likelihood ratio function  is (see \cite{LS})
\begin{align*}
V\left(\vartheta ,X^T\right)&=\exp\left\{\frac{1}{\sigma ^2}\int_{0}^{T}
f\left(\vartheta t\right)m\left(\vartheta,t \right){\rm d}X_t\right.\\
&\qquad \qquad  \left. -\frac{1}{2\sigma ^2}\int_{0}^{T}
f\left(\vartheta t\right)^2m\left(\vartheta,t \right)^2{\rm d}t
\right\},\qquad \quad \vartheta \in 
\Theta .
\end{align*}
The MLE $\hat\vartheta _T$ is defined by the relation
\begin{align*}
V\left(\hat\vartheta _T,X^T\right)=\sup_{\vartheta \in\Theta }V\left(\vartheta
,X^T\right) .
\end{align*}

We are interested in the asymptotic behavior of this estimator. The
asymptotics providing the consistency of estimators for this model of
observations can be, for example, the
following: 
\begin{description}{\it 
\item[a)] $\sigma \rightarrow 0$,   $b\rightarrow 0$ and $T$ is fixed.,
 \item[b)] $\sigma \rightarrow   0$,  $b $ and   $T$ are fixed,
\item[c)] $T\rightarrow \infty $, $\sigma $ and $b$ are fixed. }
\end{description}
Note that the asymptotic $A\rightarrow \infty $ in \eqref{0} can
be reduced to case b).
In all three cases this problem of parameter
estimation is regular and in such situations these estimators are usually
asymptotically normal with the {\it natural normalization} by the Fisher information:
\begin{align*}
\sqrt{{\rm I}_T\left(\vartheta \right)}\left( \hat\vartheta _T-\vartheta
\right) \Longrightarrow  {\cal N}\left(0,1\right),\qquad \quad \sqrt{{\rm
    I}_T\left(\vartheta \right)}\left( \tilde\vartheta _T-\vartheta 
\right) \Longrightarrow  {\cal N}\left(0,1\right) .
\end{align*}
Here  ${\rm I}_T\left(\vartheta \right)
$ is the Fisher information
\begin{align}
\label{F}
{\rm I}_T\left(\vartheta\right)=\frac{1}{\sigma ^2}\int_{0}^{T}
\left[f\left( \vartheta t\right)\dot m\left(\vartheta
  ,t\right)+t  f'\left(  \vartheta t\right)m\left(\vartheta
  ,t\right) \right]^2{\rm d}t .
\end{align}
Here and in the sequel dot means differentiation  w.r.t. $\vartheta $ and prim means
differentiation w.r.t. $t$. For example, $f'\left(\vartheta  t\right)=\left.\frac{{\rm
    d}f\left(s\right)}{{\rm d}s}\right|_{s=\vartheta t}$. 

In the work \cite{ChK1} we studied the asymptotic a) with $\sigma
=b=\varepsilon \rightarrow 0$. It was shown that the
MLE and BE are consistent, asymptotically normal 
\begin{align*}
\frac{\hat\vartheta _\varepsilon -\vartheta }{\varepsilon }\Rightarrow {\cal
  N}\left(0,{\rm I}\left(\vartheta\right)^{-1}\right),\qquad \quad
\frac{\tilde\vartheta _\varepsilon -\vartheta }{\varepsilon }\Rightarrow
     {\cal N}\left(0,{\rm I}\left(\vartheta \right)^{-1}\right)
\end{align*}
and are asymptotically efficient. 

In the present work we consider the asymptotic of the case b), i.e., we put
$\sigma =\varepsilon \rightarrow 0$ and the coefficient $b>0$ and $T>0$
 keep fixed. This model of observation has some interesting features. Let us
 see what happens with the Fisher information \eqref{F}.  The
 first strange result  is the following limit: for all $t\in(0,T]$
\begin{align*}
\lim_{\varepsilon \rightarrow 0}\left[f\left( \vartheta t\right)\dot m\left(\vartheta
  ,t\right)+t  f'\left(  \vartheta t\right)m\left(\vartheta
  ,t\right)\right]=0 .
\end{align*}
This means that 
$$
\varepsilon ^2 {\rm I}_T\left(\vartheta\right)=\int_{0}^{T} \left(\frac{\partial
}{\partial  \vartheta }\left[f\left(\vartheta t\right)m\left(\vartheta
  ,t\right)\right]    \right)^2{\rm d}t\longrightarrow  0
$$
as $\varepsilon \rightarrow 0$.
Further, we have the convergence {\it in distribution}
\begin{align*}
\varepsilon ^{-1/2}  \left[f\left( \vartheta t\right)\dot m\left(\vartheta
  ,t\right)+t  f'\left(  \vartheta t\right)m\left(\vartheta
  ,t\right)\right]\Longrightarrow \sqrt{b}\, tf'\left(\vartheta t\right)\xi _t,
\end{align*}
where $\left\{\xi _t,t\in (0,T]\right\}$ is a family of  independent  Gaussian   random
  variables, $\xi _t\sim{\cal N}\left(0,\frac{1}{2}\right)$. The  limit integral
  of the normalized Fisher information $ 
  \varepsilon {\rm I}_T\left(\vartheta\right) $ is equal (formally) to  
\begin{align*}
{b}\,\int_{0}^{T} t^2f'\left(\vartheta t\right)^2\xi _t^2{\rm d}t
\end{align*}
but this integral by the well known reason does not exist.  It is shown that   the following limit 
\begin{align*}
\varepsilon {\rm I}_T\left(\vartheta\right)\longrightarrow {\rm
  I}_0\left(\vartheta\right)= \frac{b}{2} \int_{0}^{T} t^2f'\left(\vartheta
t\right)^2{\rm d}t
\end{align*}
holds. 

The main result of the work is  the 
asymptotically normality of  the MLE:
\begin{align*}
\frac{\hat \vartheta _\varepsilon -\vartheta  }{\sqrt{\varepsilon
}}\Longrightarrow {\cal N}\left(0,{\rm I}_0\left(\vartheta\right)^{-1}
\right). 
\end{align*}
In the next section we give some auxiliary results from the Kalman filtration
and show that the limit model ($\varepsilon =0$) admits estimation of the
parameter $\vartheta $ without error.

\section{Auxiliary results}

The process $m\left(\vartheta ,t\right),0\leq t\leq T$ satisfies the
following equations of Kalman-Bucy filtration \cite{K-B61} (see details in
\cite{LS}, Theorem 10.1) 
\begin{align}
\label{kf}
{\rm d}m\left(\vartheta ,t\right)&=-am\left(\vartheta ,t\right){\rm d}t\nonumber\\
&\qquad \qquad +
\frac{\gamma \left(\vartheta ,t\right)f\left(\vartheta
  t\right)}{\varepsilon  ^2} \left[{\rm d}X_t-f\left(\vartheta
  t\right)m\left(\vartheta ,t\right){\rm d}t\right], 
\end{align}
where the function $\gamma \left(\vartheta ,t\right)=\Ex_\vartheta
\left(m\left(\vartheta ,t\right)-Y_t \right)^2 $ is solution of the Riccati
equation 
\begin{align}
\label{eR}
\frac{\partial \gamma \left(\vartheta ,t\right)}{\partial t}=
-2a\gamma \left(\vartheta ,t\right) -\frac{\gamma \left(\vartheta
  ,t\right)^2f\left(  \vartheta t\right)^2}{\varepsilon  ^2} +b^2,\qquad
\gamma \left(\vartheta ,0\right)=0. 
\end{align}

Therefore for the derivative $\dot m_t\left(\vartheta ,t\right)$ we obtain 
\begin{align}
&{\rm d}\dot m\left(\vartheta ,t\right)=-\left[a+ \frac{\gamma
    \left(\vartheta ,t\right) f\left(  \vartheta t\right)^2}{\varepsilon 
    ^2}\right]\dot m\left(\vartheta ,t\right){\rm d}t\nonumber\\ &\qquad 
\quad-h\left(\vartheta ,t\right)m\left(\vartheta ,t\right){\rm d}t +
g\left(\vartheta ,t\right)\left[{\rm d}X_t-f\left(  \vartheta
  t\right)m\left(\vartheta ,t\right){\rm d}t\right],
\label{dm}
\end{align}
where  $\dot m\left(\vartheta ,0\right)=0$ and we denoted
\begin{align*}
 g\left(\vartheta ,t\right)&=\frac{\dot \gamma
  \left(\vartheta ,t\right)f\left(  \vartheta t\right)+ t\gamma
  \left(\vartheta ,t\right) f'\left( \vartheta t\right)}{\varepsilon ^2 },\\
h\left(\vartheta ,t\right)&=\frac{t\gamma \left(\vartheta ,t\right)f\left(  \vartheta
  t\right) f'\left(  \vartheta
  t\right)}{\varepsilon  ^2} .
\end{align*}
For the derivative  $\dot \gamma \left(\vartheta ,t\right) $ we have the 
equation
\begin{align}
\frac{\partial \dot \gamma \left(\vartheta ,t\right)}{\partial t }&=
-2\left[a +\frac{\gamma \left(\vartheta
  ,t\right)f\left(  \vartheta t\right)^2}{\varepsilon  ^2}\right]\dot \gamma
\left(\vartheta ,t\right)\nonumber\\
&\qquad  -\frac{2 t \gamma \left(\vartheta
  ,t\right)^2f\left( \vartheta t\right) f'\left( \vartheta
  t\right)}{\varepsilon  ^2} ,\quad \dot \gamma \left(\vartheta ,0\right)=0. 
\label{gm}
\end{align}

These equations are obtained by the formal differentiation  but this
derivation can be justified by the standard methods. The 
both functions (Gaussian $m\left(\vartheta ,t\right)$ and deterministic
$g\left(\vartheta ,t\right)$) are infinitely differentiable.

The equations \eqref{dm} and \eqref{gm} are linear and their solutions can be
written explicitly. Let us denote 
\begin{align*}
q\left(\vartheta ,t\right)=a+ \frac{\gamma_*
    \left(\vartheta ,t\right) f\left(  \vartheta
    t\right)^2}{\varepsilon},\qquad \gamma_* 
    \left(\vartheta ,t\right)=\frac{\gamma
    \left(\vartheta ,t\right)}{\varepsilon }.
\end{align*}
Then we have
\begin{align*}
\dot m\left(\vartheta ,t\right)&=-\int_{0}^{t}e^{-\int_{s}^{t}q\left(\vartheta
  ,v\right){\rm d}v } h\left(\vartheta ,s\right)m\left(\vartheta ,s\right){\rm
  d}s\\
&\qquad  +\int_{0}^{t}e^{-\int_{s}^{t}q\left(\vartheta
  ,v\right){\rm d}v }
g\left(\vartheta ,s\right)\left[{\rm d}X_s-f\left(  \vartheta
  s\right)m\left(\vartheta ,s\right){\rm d}s\right]
\end{align*}
and
\begin{align}
\label{sgm}
\dot \gamma \left(\vartheta ,t\right)&=-2\int_{0}^{t}se^{-2\int_{s}^{t}q\left(\vartheta 
  ,v\right){\rm d}v }   \gamma_* \left(\vartheta
  ,s\right)^2 f\left( \vartheta s\right)f'\left( \vartheta s\right){\rm
  d}s.
\end{align}

Let us see how can be constructed a consistent estimator of $\vartheta $ in the
case of asymptotic $\varepsilon \rightarrow 0 $ by the  observations
$X^T$. The model of observations is 
\begin{align*}
{\rm d} X_t&=f\left(  \vartheta t\right)Y_t{\rm d}t+\varepsilon {\rm
  d}W_t,\quad X_0=0,\quad 0\leq t\leq T,\\ 
{\rm d}Y_t&=-aY_t{\rm d}t+b{\rm d}V_t,\qquad\quad Y_0=y_0.
\end{align*}

Suppose that $\varepsilon =0$ (limit system) and construct an estimator of $\vartheta $
without error. Hence
\begin{align}
\label{lm}
\frac{{\rm d} x_t}{{\rm d}t}&=f\left(  \vartheta t\right)Y_t,\qquad x_0=0,\\
{\rm d}Y_t&=-aY_t{\rm d}t+b{\rm d}V_t,\quad Y_0=y_0 \nonumber
\end{align}
and we have to estimate $\vartheta $ by the observations $x^T=\left(x_t,0\leq
t\leq T\right)$. 
 Let us put $z_t= \frac{{\rm d} x_t}{{\rm d}t}$, then
\begin{align*}
{\rm d}z_t&= \vartheta   f'\left(  \vartheta t\right)Y_t{\rm
  d}t+f\left(  \vartheta t\right){\rm d}Y_t\\
&=\left[  \vartheta
 f'\left(  \vartheta t\right)-af\left(  \vartheta
t\right)\right]Y_t{\rm d}t+bf\left(  \vartheta t\right){\rm d}V_t. 
\end{align*}
By the It\^o formula
\begin{align*}
z_t^2=2\int_{0}^{t}z_s{\rm d}z_s+{b^2}{}\int_{0}^{t} f\left( 
\vartheta s\right)^2{\rm d}s.
\end{align*}
Hence the function
\begin{align*}
\Psi \left(t\right)=z_t^2-2\int_{0}^{t}z_s{\rm
  d}z_s=b^2\int_{0}^{t} f\left( 
\vartheta s\right)^2{\rm d}s
\end{align*}
is deterministic and for any $t\in (0,T]$ the observed $\Psi \left(t\right)$
  defines $\vartheta $ without error. This means that if we have the limit
  model \eqref{lm}, then the measures corresponding to the observations are
  singular. 

Suppose that $f\left(\vartheta t\right)=A\cos\left(\vartheta t\right)$. Then
\begin{align*}
\int_{0}^{t} f\left( 
\vartheta s\right)^2{\rm d}s=\frac{A^2t}{2}+\frac{A^2}{4\vartheta
}\sin\left(2\vartheta t\right) 
\end{align*}
and
\begin{align*}
\tilde \Psi \left(t\right)=\frac{\Psi
  \left(t\right)}{b^{2}}-\frac{A^2t}{2}=\frac{A^2}{4\vartheta
}\sin\left(2\vartheta t\right). 
\end{align*}
If we denote
\begin{align*}
\tau =\arg\inf_{t>t_0}\tilde \Psi \left(t\right)=0
\end{align*}
then $\vartheta =\frac{\pi }{2\tau }$.

Therefore if $\varepsilon \rightarrow 0$ then the consistent
  estimation is possible. Of course, we cannot differentiate the observations
  $X^T$ w.r.t. $t$ but we can do it ``asymptotically'' with the help of the
  kernel. For example, let us define
\begin{align*}
\hat z_t=\frac{1}{\varphi _\varepsilon }\int_{0}^{T}K\left(\frac{s-t}{\varphi
  _\varepsilon }\right){\rm d}X_s. 
\end{align*}
Here the kernel $K\left(\cdot \right)$ satisfies the usual conditions:
\begin{align*}
K\left(u\right)\geq 0,\qquad \int_{c_1}^{c_2}K\left(u\right){\rm d}u=1
\end{align*}
and $K\left(u\right)=0$ for $u\not \in \left(c_1,c_2\right)$. Moreover we
suppose that $K\left(u\right)$ is continuously differentiable. 

Then for $t\in \left(0,T\right)$ and small $\varepsilon $ we have 
\begin{align*}
z_t=f\left(\vartheta t\right)Y_t= \frac{1 }{\varphi _\varepsilon
}\int_{0}^{T}K\left(\frac{s-t}{\varphi _\varepsilon }\right)f\left(\vartheta t\right)Y_t{\rm
  d}s.
\end{align*}
Hence we can write 
\begin{align*}
\Ex_\vartheta \left(\hat z_t- z_t\right)^2&=\Ex_\vartheta \left(
\frac{\varepsilon }{\varphi _\varepsilon
}\int_{0}^{T}K\left(\frac{s-t}{\varphi _\varepsilon }\right){\rm
  d}W_s\right)^2\\
&\qquad +\Ex_\vartheta \left(
\frac{1}{ \varphi _\varepsilon
}\int_{0}^{T}K\left(\frac{s-t}{\varphi _\varepsilon
}\right)\left[f\left(\vartheta s\right)Y_s-f\left(\vartheta t\right)Y_t
  \right]{\rm d}s\right)^2.
\end{align*}
Below we put $s=t+\varphi _\varepsilon u$
\begin{align*}
\Ex_\vartheta \left(
\frac{\varepsilon }{\varphi _\varepsilon
}\int_{0}^{T}K\left(\frac{s-t}{\varphi _\varepsilon }\right){\rm
  d}W_s\right)^2=\frac{\varepsilon ^2}{\varphi _\varepsilon }\int_{-\frac{t}{\varphi
    _\varepsilon }}^{\frac{T-t}{\varphi _\varepsilon }} K\left(u\right)^2{\rm d}u
\end{align*}
and
\begin{align*}
&\frac{1}{ \varphi _\varepsilon
}\int_{0}^{T}K\left(\frac{s-t}{\varphi _\varepsilon
}\right)\left[f\left(\vartheta s\right)Y_s-f\left(\vartheta t\right)Y_t
  \right]{\rm d}s\\
&\qquad \qquad =\frac{1}{ \varphi _\varepsilon
}\int_{0}^{T}K\left(\frac{s-t}{\varphi _\varepsilon
}\right)\left[f\left(\vartheta s\right)-f\left(\vartheta t\right)
  \right]Y_s{\rm d}s\\ 
& \qquad\qquad  \quad +\frac{f\left(\vartheta t\right)}{ \varphi _\varepsilon
}\int_{0}^{T}K\left(\frac{s-t}{\varphi _\varepsilon
}\right)\left[Y_s-Y_t  \right]{\rm d}s\\
&\qquad \qquad =\int_{-\frac{t}{\varphi
    _\varepsilon } }^{\frac{T-t}{\varphi _\varepsilon
  }}K\left(u\right)\left[f\left(\vartheta \left(t+\varphi _\varepsilon u
    \right)\right)-f\left(\vartheta t\right) 
  \right]Y_{t+\varphi _\varepsilon u}{\rm d}u\\
&\qquad\qquad  \quad +{f\left(\vartheta t\right)}\int_{-\frac{t}{\varphi
    _\varepsilon }}^{\frac{T-t}{\varphi _\varepsilon
  }}K\left(u\right)\left[Y_{t+\varphi _\varepsilon u}-Y_t  \right]{\rm d}s\\ 
&\qquad \qquad =\varphi _\varepsilon \vartheta \int_{-\frac{t}{\varphi 
    _\varepsilon } }^{\frac{T-t}{\varphi _\varepsilon
  }}uK\left(u\right)f'\left(\vartheta \tilde t _\varepsilon
  \right)Y_{t+\varphi _\varepsilon u}{\rm d}u\\ 
&\qquad\qquad  \quad +\sqrt{\varphi _\varepsilon }{f\left(\vartheta
    t\right)}\int_{-\frac{t}{\varphi 
    _\varepsilon }}^{\frac{T-t}{\varphi _\varepsilon
  }}K\left(u\right)\left[\frac{Y_{t+\varphi _\varepsilon u}-Y_t
    }{\sqrt{\varphi _\varepsilon }} \right]{\rm d}s. 
\end{align*}
For the process $Y_t$ we have 
\begin{align*}
\frac{Y_{t+\varphi _\varepsilon u}-Y_t }{\sqrt{\varphi _\varepsilon
}}=\frac{-a}{\sqrt{\varphi _\varepsilon }}\int_{ t}^{t+\varphi _\varepsilon u}
Y_s{\rm d}s+b\frac{V_{t+\varphi _\varepsilon u}-V_t}{\sqrt{\varphi _\varepsilon }} .
\end{align*}
Hence
\begin{align*}
\Ex_{\vartheta }\left( \frac{Y_{t+\varphi _\varepsilon u}-Y_t }{\sqrt{\varphi _\varepsilon
}} \right)^2\leq C
\end{align*}
with some constant $C=C\left(u,t\right)>0$. 
Therefore we obtain the following estimate for the error
\begin{align*}
\Ex_\vartheta \left(\hat z_t- z_t\right)^2\leq C_1\frac{\varepsilon
  ^2}{\varphi _\varepsilon }+C_2 \varphi _\varepsilon \leq C\varepsilon 
\end{align*}
if we take the {\it optimal} choice  $\varphi _\varepsilon  =\varepsilon
$. This means that
\begin{align*}
\hat z_t=f\left(\vartheta t\right)Y_t+O\left(\sqrt{\varepsilon }\right).
\end{align*}
The function $\Psi \left(t\right)$ can be estimated as follows
\begin{align*}
\hat \Psi _\varepsilon \left(t\right)=\hat z_t^2-2\sum_{k=0}^{K-1}\hat
z_{t_k}\left[ \hat z_{t_{k+1}  }-\hat z_{t_k}\right] .
\end{align*}
Here $t_k=\frac{kt}{K}$ and it can be shown that for a
special choice $K=K\left(\varepsilon \right)\rightarrow \infty $  we obtain
$\hat \Psi _\varepsilon \left(t\right)\rightarrow \Psi _\varepsilon
\left(t\right) $. Using the standard arguments we can verify the consistency
of   the estimator $\vartheta _\varepsilon ^*$ defined by the
relation
\begin{align*}
\hat \Psi _\varepsilon \left(t\right)=b^2\int_{0}^{t}f\left(\vartheta
_\varepsilon^* s\right)^2{\rm d}s .
\end{align*}
Of course, this estimator has a bad rate of convergence.

\section{Main result}

The MLE and BE we denote as $\hat\vartheta
_\varepsilon $ and $\tilde\vartheta_\varepsilon $  respectively. We have the model of observations
\begin{align*}
{\rm d}X_t&=f\left(  \vartheta t\right)Y_t{\rm d}t+\varepsilon  {\rm
  d}W_t,\qquad X_0=0,\quad 0\leq t\leq T, \\ 
{\rm d}Y_t&=-aY_t{\rm d}t+b {\rm
  d}V_t,\qquad\qquad\qquad Y_0=y_0.
\end{align*}
Our goal is to estimate $\vartheta $ and to describe the properties of the
estimators as $\varepsilon \rightarrow 0$.   We suppose
that the periodic function $f\left(t\right)$ is positive and two times continuously
differentiable.  We denote $f'\left(t\right)$ the derivative and put 
\begin{align}
\label{e} 
\kappa=\inf_{0\leq t\leq 1}f\left(t\right)>0,\qquad K=\sup_{0\leq
  t\leq 1} f\left(t\right)<\infty .
\end{align}
 This allows us to avoid
the situation, where $Y_t$ is multiplied by 0. Recall that we suppose that 
 $a,b>0, y_0$ are known, the conditions \eqref{e} are fulfilled  and  the Fisher
information is
\begin{align*}
{\rm I}_0\left(\vartheta \right)=\frac{b}{2}\int_{0}^{T}t^2f'\left(\vartheta t\right)^2{\rm d}t.
\end{align*}
The main result of this work is the following theorem.
\begin{theorem}
\label{T1} The MLE $\hat\vartheta_\varepsilon  $ is consistent and
asymptotically normal
\begin{align}
\label{an}
\frac{\hat\vartheta_\varepsilon -\vartheta }{\sqrt{\varepsilon
}}\Longrightarrow {\cal N}\left(0,{\rm I}_0\left(\vartheta \right)^{-1}\right).
\end{align}
\end{theorem}
{\bf Proof.}  Let us denote $M\left(\vartheta ,t\right)=f\left(\vartheta
t\right)m\left(\vartheta ,t\right)$, where $m\left(\vartheta ,t\right) $ is
solution of the equation \eqref{kf}. We have to study the log-likelihood ratio
\begin{align*}
\ln V\left(\vartheta ,X^T\right)=\frac{1}{\varepsilon ^2}\int_{0}^{T}
M\left(\vartheta ,t\right){\rm d}X_t-\frac{1}{2\varepsilon ^2} \int_{0}^{T}
M\left(\vartheta ,t\right)^2{\rm d}t,\; \vartheta \in\Theta =\left(\alpha
,\beta \right). 
\end{align*} 
We have the relation
\begin{align*}
\dot M\left(\vartheta ,t\right)=tf'\left(\vartheta t\right)m\left(\vartheta
,t\right)+f\left(\vartheta t\right)\dot m\left(\vartheta ,t\right)
\end{align*}
and we need to know the asymptotics of the random processes $m\left(\vartheta
,t\right),0\leq t\leq T$ and $\dot m\left(\vartheta ,t\right),0\leq t\leq T$ as $\varepsilon \rightarrow
0$. 
 Introduce the
function $\gamma _*\left(\vartheta ,t\right)=\varepsilon ^{-1}\gamma
\left(\vartheta ,t\right)$ and note that 
\begin{align*}
{\rm d}m\left(\vartheta ,t\right)=-am\left(\vartheta ,t\right){\rm d}t  +
{\gamma_* \left(\vartheta ,t\right)f\left(  \vartheta
  t\right)}{\rm d}\bar W_t.
\end{align*}
Here $\bar W_t$ is the {\it innovation} Wiener process 
\begin{align*}
\bar W_t=\frac{1}{\varepsilon }\left[{\rm d}X_t-M\left(\vartheta ,t\right){\rm
    d}t\right]. 
\end{align*}
The asymptotics of the solution of Riccati equation is given by the following
lemma.
\begin{lemma}
\label{L1} 
For any $t_0\in (0,T]$ we have the convergence
\begin{align}
\label{asym}
\sup_{t_0\leq t\leq T}\left|\gamma \left(\vartheta
,t\right)-\frac{b\,\varepsilon }{f\left(\vartheta t\right)}\right| =O\left(\varepsilon ^2\right).
\end{align}
\end{lemma}
{\bf Proof.} Recall that the function 
 $\gamma \left(\vartheta ,t\right)=\Ex_{\vartheta }\left(m\left(\vartheta
,t\right)-Y_t\right)^2$
satisfies the Riccati equation \eqref{eR}
\begin{align*}
\frac{\partial \gamma \left(\vartheta ,t\right)}{\partial t }=
-2a\gamma\left(\vartheta ,t\right) -{\gamma_* \left(\vartheta
  ,t\right)^2f\left(\vartheta t\right)^2} +b^2,\qquad \gamma \left(\vartheta
,0\right)=0. 
\end{align*}
To verify the convergence \eqref{asym} we introduce the equation
\begin{align*}
\frac{\partial \gamma ^*\left(t\right)}{\partial t }=
-2a\gamma^*\left(t\right) -\frac{{\gamma^* \left(t\right)^2\kappa
    ^2}}{\varepsilon ^2} +b^2,\qquad \gamma^* \left(0\right)=0.
\end{align*}
and note that by the comparison theorem for ordinary differential equations we
have the relation
\begin{align*}
\gamma \left(\vartheta,t \right)\leq \gamma^* \left(t\right),\qquad 0\leq t\leq T.
\end{align*}
The solution $ \gamma^* \left(t\right)$ can be written explicitly \cite{A83}
\begin{align*}
\gamma^* \left(t\right)=e^{-2rt}\left[\frac{1}{\gamma^*\left(0\right)-\hat\gamma
  }+\frac{\kappa ^2}{2r\varepsilon ^2} \left(1-e^{-2rt}\right)\right]^{-1}+\hat\gamma .
\end{align*}
Here we denoted 
\begin{align*}
r=\left(a^2+\frac{b^2\kappa ^2}{\varepsilon ^2}\right)^{1/2},\qquad \quad 
\hat\gamma =\frac{a\varepsilon ^2}{\kappa ^2}\left(\sqrt{1+\frac{b^2\kappa
    ^2}{a^2\varepsilon ^2}} -1\right).
\end{align*}
It is easy to see that for any $t_0\in(0,T]$ we have the representations
\begin{align*}
&r=\frac{b\kappa }{\varepsilon }\left(1+O\left(\varepsilon \right)\right),\qquad
\hat\gamma =\frac{b\varepsilon }{\kappa }\left(1+O\left(\varepsilon \right)\right),\\
&\sup_{t_0\leq t\leq T}\left|\gamma^* \left(t\right)-\frac{b\varepsilon
}{\kappa }\right|=O\left(\varepsilon^2 \right).
\end{align*}
Hence for $t>t_0$ and $0<\varepsilon \leq \varepsilon _0$ with some
$\varepsilon _0>0$ we have 
\begin{align*}
0\leq \gamma_* \left(\vartheta ,t\right)=\frac{\gamma \left(\vartheta ,t\right)}{\varepsilon }\leq \frac{2b
}{\kappa }.
\end{align*}
Using the similar arguments   we obtain the following estimate from below
\begin{align*}
\gamma_* \left(\vartheta ,t\right) \geq \frac{b }{2K}.
\end{align*}
We have the relation
\begin{align*}
\gamma \left(\vartheta ,t\right)-\gamma \left(\vartheta ,t_0\right)+a\int_{t_0}^{t} \gamma \left(\vartheta
,s\right){\rm d}s=-\int_{t_0}^{t} \gamma_* \left(\vartheta
,s\right)^2f\left(\vartheta s\right)^2{\rm d}s +b^2\left(t-t_0\right),
\end{align*}
where the left hand part tends to zero. 

Hence, we verified \eqref{asym} and can write 
\begin{align*}
{\rm d}m\left(\vartheta ,t\right)=-am\left(\vartheta ,t\right){\rm d}t  +
b\left(1+o\left(1\right)\right){\rm d}\bar W_t,\qquad m\left(\vartheta
,0\right)=y_0. 
\end{align*}

Below we will use several time the following technical elementary lemma.
\begin{lemma}
\label{L0} Suppose that the functions $F\left(t\right), t\in [0,T]$ and
  $G\left(t\right), t\in\left[0,T\right]$  are continuously differentiable,
  the function $F(0)=0$, $F\left(t\right)>0, t\in (0,T]$ and
  $\varepsilon \rightarrow 0$, then we have the estimate
\begin{align*}
N_\varepsilon \left(t\right)=\int_{0}^{t}e^{-\frac{1}{\varepsilon }\int_{s}^{t}F\left(v\right){\rm
    d}v}G\left(s\right){\rm d}s=\varepsilon\,\frac{
  G\left(t\right)}{F\left(t\right)}\left(1+O\left(\varepsilon \right) \right) 
\end{align*}
for any $t>0$. 
\end{lemma}
{\bf Proof.} Let us take some  (small) $t_0>0$ and  $t_1\in
\left(t_0,t\right)$  and denote 
\begin{align*}
\inf_{t_0\leq s\leq T}F\left(s\right)=c_1>0,\qquad \sup_{0\leq s\leq
  T}F\left(s\right)=C_1<\infty ,\quad \sup_{0\leq s\leq
  T}\left|G\left(s\right)\right|=C_2<\infty . 
\end{align*}
We have the estimate
\begin{align*}
J_0=\int_{0}^{t_0}e^{-\frac{1}{\varepsilon }\int_{s}^{t}F\left(v\right){\rm
    d}v}\left|G\left(s\right)\right|{\rm d}s \leq C_2t_0
e^{-\frac{1}{\varepsilon }\int_{t_0}^{t}F\left(v\right){\rm     d}v} \leq C_2t_0
e^{-\frac{\left(t-t_0\right)}{\varepsilon }}.
\end{align*}
Then we can write
\begin{align*}
J_1=\int_{t_0}^{t_1}e^{-\frac{1}{\varepsilon }\int_{s}^{t}F\left(v\right){\rm
    d}v}\left|G\left(s\right)\right|{\rm d}s \leq
\int_{t_0}^{t_1}e^{-\frac{\left(t-s\right)c_1}{\varepsilon } }
\left|G\left(s\right)\right|{\rm d}s
\leq  \frac{C_2\varepsilon }{c_1}\;e^{-\frac{\left(t-t_1\right)c_1}{\varepsilon }}.
\end{align*}
Below we change the variables $s=t-u\varepsilon $,  $v=t-q\varepsilon $ and
use the Taylor expansion 
\begin{align*}
J_2&=\int_{t_1}^{t}e^{-\frac{1}{\varepsilon }\int_{s}^{t}F\left(v\right){\rm
    d}v}G\left(s\right){\rm d}s=-\varepsilon \int_{\frac{t-t_1}{\varepsilon
}}^{0} e^{-\frac{1}{\varepsilon }\int_{t-u\varepsilon }^{t}F\left(v\right){\rm
    d}v} G\left(t-u\varepsilon \right){\rm d}u\\
&= \varepsilon G\left(t \right)\int_{0}^{\frac{t-t_1}{\varepsilon
}} \exp\left\{ -\frac{1}{\varepsilon } \int_{t-u\varepsilon
}^{t}\left[F\left(t\right)+\left(t-v\right)F'\left(\tilde t\right)\right] {\rm
    d}v  \right\}\left(1+O\left(\varepsilon \right)\right)\\
&= \varepsilon G\left(t \right)\int_{0}^{\frac{t-t_1}{\varepsilon
}} \exp\left\{ -u F\left(t\right)+O\left(\varepsilon \right)     \right\}{\rm
  d}u \left(1+O\left(\varepsilon \right)\right)= \varepsilon \frac{G\left(t \right)}{F\left(t\right)} \left(1+O\left(\varepsilon \right)\right)
\end{align*}

The asymptotics of $\dot M\left(\vartheta ,t\right)$ is described in the next
lemma. 
\begin{lemma}
\label{L2} For any $t\in (t_0,T]$ we have the limits
\begin{align}
\label{l2}
m\left(\vartheta ,t\right)\rightarrow Y_t,\qquad \dot m\left(\vartheta
,t\right)\longrightarrow - \frac{t f'\left(\vartheta 
  t\right)}{f\left(\vartheta t\right)}Y_t 
\end{align}
as $\varepsilon \rightarrow 0$ and therefore $\dot M\left(\vartheta ,t\right)\rightarrow 0 $.
\end{lemma}
{\bf Proof.} The first convergence follows immediately from $ \Ex_{\vartheta
}\left(m\left(\vartheta ,t\right)-Y_t\right)^2=\varepsilon \gamma
_*\left(\vartheta ,t\right)\rightarrow 0$, i.e.,  we have the mean square convergence
$m\left(\vartheta ,t\right) \rightarrow Y_t$ uniformly on $t\in
\left[t_0,T\right]$ for any $t_0\in (0,T]$.    The
  derivative $\dot m\left(\vartheta ,t\right)$   satisfies the equation
\begin{align*}
\dot m\left(\vartheta ,t\right)&=-\frac{1}{\varepsilon
}\int_{0}^{t}e^{-\int_{s}^{t}q\left(\vartheta 
  ,v\right){\rm d}v } h_*\left(\vartheta ,s\right)m\left(\vartheta ,s\right){\rm
  d}s\\
&\qquad  +\int_{0}^{t}e^{-\int_{s}^{t}q\left(\vartheta
  ,v\right){\rm d}v }
g_*\left(\vartheta ,s\right){\rm d}\bar W_s,
\end{align*}
where we denoted $g_*\left(\vartheta ,t\right)={\dot \gamma_*
  \left(\vartheta ,t\right)f\left(  \vartheta t\right)+ t\gamma_*
  \left(\vartheta ,t\right) f'\left( \vartheta t\right)} $ and
$h_*\left(\vartheta ,t\right)={t\gamma_* \left(\vartheta ,t\right)f\left( 
  \vartheta 
  t\right) f'\left(  \vartheta
  t\right)}$. Here $\dot \gamma_*
  \left(\vartheta ,t\right)=\varepsilon ^{-1} \dot \gamma
  \left(\vartheta ,t\right)$.
Note that for the values $v\in \left[s,t\right]$ with $s>t_0$ and
$\left|t-s\right|\leq C\varepsilon $ we have 
\begin{align*}
q\left(\vartheta ,v\right)=\frac{1}{\varepsilon }\left[a\varepsilon +\gamma _*\left(\vartheta
,v\right)f\left(\vartheta v\right)^2 \right] =\frac{bf\left(\vartheta
  t\right)\left(1+o\left(1\right)\right)}{\varepsilon } 
\end{align*}
The derivative $\dot\gamma_* \left(\vartheta ,t\right) $ according to  Lemma
\ref{L0} and the 
equation \eqref{sgm} has the following asymptotics (below $f_t=f\left(\vartheta _0t\right)$)
\begin{align*}
\dot\gamma \left(\vartheta
,t\right)&=-2\int_{0}^{t}s\exp\left\{-2\frac{\left(t-s\right){bf_t}}{\varepsilon}\right\}
\gamma _*\left(\vartheta s\right)^2f\left(\vartheta s\right)f'\left(\vartheta
s\right){\rm d}s \left(1+O\left(\varepsilon
\right)\right)\\ &=-\frac{2\varepsilon }{bf\left(\vartheta
  t\right)}\int_{0}^{\frac{tbf_t}{\varepsilon }}e^{-2u}\gamma
_*\left(\vartheta,s_u\right)^2f\left(\vartheta s_u\right)f'\left(\vartheta s_u
\right) {\rm d}u \left(1+O\left(\varepsilon
\right)\right)\\ &=-\frac{2\varepsilon \gamma
  _*\left(\vartheta,t\right)^2f'\left(\vartheta t
  \right)}{b}\int_{0}^{\frac{tbf_t}{\varepsilon }}e^{-2u} {\rm d}u
\left(1+O\left(\varepsilon \right)\right)\\ &=-\frac{\varepsilon
  bf'\left(\vartheta t \right)}{f\left(\vartheta t \right)}
\left(1+O\left(\varepsilon \right)\right),
\end{align*}
where we put $s=s_u=t-\frac{u\varepsilon }{bf_t}$ and used the Taylor
formula. Further
\begin{align*}
g_*\left(\vartheta ,t\right)&=-\frac{{tbf'\left(\vartheta
    t\right)}}{f\left(\vartheta t\right) } \left(1+O\left(\varepsilon
  \right)\right)+\frac{{tbf'\left(\vartheta
    t\right)}}{f\left(\vartheta t\right) } \left(1+O\left(\varepsilon
  \right)\right)=O\left(\varepsilon
  \right).
\end{align*}
This allows us to write
\begin{align*}
\Ex_{\vartheta _0}\left(\int_{0}^{t}e^{-\int_{s}^{t}q\left(\vartheta
  ,v\right){\rm d}v } g_*\left(\vartheta ,s\right){\rm d}\bar
W_s\right)^2=O\left(\varepsilon^2 \right)
\end{align*}
and
\begin{align*}
\dot m\left(\vartheta ,t\right)&=-\frac{1}{\varepsilon
}\int_{0}^{t}e^{-\int_{s}^{t}q\left(\vartheta ,v\right){\rm d}v } {s\gamma_*
  \left(\vartheta ,s\right)f\left( \vartheta s\right) f'\left( \vartheta
  s\right)}m\left(\vartheta ,s\right){\rm d}s+O\left(\varepsilon
\right)\\ 
&=-\frac{t\gamma_* \left(\vartheta ,t\right)f\left( \vartheta
  t\right) f'\left( \vartheta t\right)}{\varepsilon
}\int_{0}^{t}e^{-\int_{s}^{t}q\left(\vartheta ,v\right){\rm d}v }
m\left(\vartheta ,s\right){\rm d}s+O\left(\varepsilon \right)\\ 
&=-\frac{t\gamma_* \left(\vartheta ,t\right) f'\left( \vartheta t\right)Y_t}{b
}+O\left(\sqrt{\varepsilon }\right)=-\frac{tb f'\left( \vartheta
  t\right)Y_t}{f \left(\vartheta t\right)
}+O\left(\sqrt{\varepsilon }\right).
\end{align*}
Here we used the relations
\begin{align*}
m\left(\vartheta ,s\right)=Y_s+O\left(\sqrt{\varepsilon
}\right)=Y_t+O\left(\sqrt{\varepsilon }\right) 
\end{align*}
which can be easily verified. 

Therefore  we have the limits \eqref{l2} and 
\begin{align*}
t f'\left(\vartheta t\right)m\left(\vartheta ,t\right)&\longrightarrow  t
f'\left(\vartheta t\right)Y_t,\qquad  f\left(\vartheta
t\right)\dot m\left(\vartheta ,t\right)\longrightarrow   - t
f'\left(\vartheta t\right)Y_t,\\
\dot M\left(\vartheta ,t\right)&\longrightarrow 0.
\end{align*}
 Hence for  the Fisher information   we obtain the   limit
\begin{align*}
{\varepsilon ^2}{\rm I}_\varepsilon \left(\vartheta \right)&=\int_{0}^{T}
\left[tf'\left(\vartheta t\right)m\left(\vartheta,t\right)+
  f\left(\vartheta t\right) \dot m\left(\vartheta ,t\right) \right] ^2{\rm 
  d}t\\ &\longrightarrow\int_{0}^{T} \left[tf'\left(\vartheta t\right)Y_t- t
  f'\left(\vartheta t\right) Y_t \right] ^2{\rm d}t
= 0.
\end{align*}

This means that we have to study the limits of the random processes
\begin{align*}
r_{t,\varepsilon }={m\left(\vartheta ,t\right)-Y_t}{},\qquad k_{t,\varepsilon }=
{\dot 
m\left(\vartheta,t\right) + \frac{{t f'\left(\vartheta
  t\right)}}{ f\left(\vartheta t\right)   }{}}Y_t.
\end{align*}
Introduce the random processes 
\begin{align*}
\zeta _{t,\varepsilon }=\int_{0}^{\frac{tbf_t}{\varepsilon }}e^{-u} {\rm
  d}W_{t,\varepsilon }\left(u\right),\qquad \xi _{t,\varepsilon
}=\int_{0}^{\frac{tbf_t}{\varepsilon }}e^{-u} {\rm d}V_{t,\varepsilon
}\left(u\right)
\end{align*}
with the  independent Wiener processes
\begin{align*}
W_{t,\varepsilon }\left(u\right)= \sqrt{\frac{bf\left(\vartheta
    t\right)}{\varepsilon }} \left(W_{t-\frac{u\varepsilon }{bf_t}} -W_t
\right),\qquad V_{t,\varepsilon }\left(u\right)= \sqrt{\frac{bf\left(\vartheta
    t\right)}{\varepsilon }} \left(V_{t-\frac{u\varepsilon }{bf_t}} -V_t
\right).
\end{align*}
For example, we have   $\Ex V_{t,\varepsilon }\left(u\right)=0$ and
$\Ex V_{t,\varepsilon }\left(u_1\right)V_{t,\varepsilon
}\left(u_2\right)=u_1\wedge u_2$.

\begin{lemma}
\label{L3} 
We have the representations
\begin{align}
\label{rep1}
r_{t,\varepsilon }&=\left(\frac{b\,\varepsilon }{f\left(\vartheta
  t\right)}\right)^{1/2}\left[\zeta _{t,\varepsilon }-\xi  _{t,\varepsilon }
  \right]\left(1+o\left(1\right)\right) ,\\
\label{rep2}
k_{t,\varepsilon }&=-\frac{tf'\left(\vartheta t\right)}{f\left(\vartheta
  t\right)}\,r_{t,\varepsilon }\left(1+o\left(1\right)\right) -{tf'\left(\vartheta
  t\right)}\sqrt{\frac{b\,\varepsilon }{f\left(\vartheta 
  t\right)^3}} \,\xi _{t,\varepsilon } \left(1+o\left(1\right)\right).
\end{align}
\end{lemma}
{\bf Proof.} 
We have 
\begin{align*}
{\rm d}r_{t,\varepsilon }=-\left[a+\frac{\gamma
    _*\left(\vartheta ,t\right)f\left(\vartheta t\right)^2 }{\varepsilon
  }\right]r_{t,\varepsilon }{\rm d}t +\gamma 
    _*\left(\vartheta ,t\right)f\left(\vartheta t\right){\rm d} W_t-b{\rm
  d}V_t,\;r_{0,\varepsilon }=0, 
\end{align*}
and
\begin{align*}
r_{t,\varepsilon }=\int_{0}^{t}e^{-\int_{s}^{t}q\left(\vartheta 
  ,v\right){\rm d}v } \left[\gamma 
    _*\left(\vartheta ,s\right)f\left(\vartheta s\right){\rm d} W_s -b{\rm
    d}V_s\right]. 
\end{align*}
Note that $\Ex_{\vartheta}r_{t,\varepsilon }=0$ and
\begin{align*}
\Ex_{\vartheta }r_{t,\varepsilon }^2=\int_{0}^{t}e^{-2\int_{s}^{t}q\left(\vartheta
  ,v\right){\rm d}v } \left[\gamma _*\left(\vartheta
  ,s\right)^2f\left(\vartheta s\right)^2+b^2\right]{\rm
  d}s=\frac{b\,\varepsilon }{f\left(\vartheta
  t\right)}\left(1+O\left(\varepsilon \right)\right).
\end{align*}
This process for $t>t_0>0$ has the following asymptotics
\begin{align*}
r_{t,\varepsilon }&=b\int_{0}^{t}e^{-\left(t-s\right)\frac{bf_t}{\varepsilon } }{\rm
  d}W_s\left(1+o\left(1\right)\right)-b\int_{0}^{t}e^{-\left(t-s\right)
\frac{bf_t}{\varepsilon } }{\rm d}V_s\left(1+o\left(1\right)\right)\\
&=\sqrt{\frac{b\,\varepsilon }{f\left(\vartheta t\right)}}\;\left[\zeta _{t,\varepsilon }-\xi _{t,\varepsilon }\right]\left(1+o\left(1\right)\right),
\end{align*}
where we changed the variables $s=t-\frac{u\varepsilon }{bf_t} $. This  proves
the first relation \eqref{rep1}.

For $k_{t,\varepsilon } $  we can write
\begin{align*}
k_{t,\varepsilon } &=-\frac{1}{\varepsilon
}\int_{0}^{t}e^{-\int_{s}^{t}q\left(\vartheta 
  ,v\right){\rm d}v } h_*\left(\vartheta ,s\right)m\left(\vartheta ,s\right){\rm
  d}s+\frac{t
  f'\left(\vartheta _0t\right) Y_t}{f\left(\vartheta t\right)
}+O\left(\varepsilon \right).
\end{align*}
Recall the estimates
\begin{align*}
\dot \gamma_*
  \left(\vartheta ,t\right)&=-\frac{tbf'\left(\vartheta
    t\right)}{f\left(\vartheta t\right)^2} \left(1+O\left(\varepsilon
  \right)\right),\qquad \gamma_* \left(\vartheta
  ,t\right)=\frac{b}{f\left(\vartheta t\right)}\left(1+O\left(\varepsilon 
  \right)\right),\\
h_*\left(\vartheta ,t\right)&={t\gamma_* \left(\vartheta ,t\right)f\left(
  \vartheta t\right) f'\left( \vartheta t\right)}=tbf'\left( \vartheta
t\right)\left(1+O\left(\varepsilon \right)\right).
\end{align*}

This allow us to write 
\begin{align*}
k_t &=-\frac{tb f'\left(\vartheta
 t\right)}{\varepsilon}\int_{0}^{t}e^{-\int_{s}^{t}q\left(\vartheta
  ,v\right){\rm d}v } m\left(\vartheta ,s\right){\rm d}s+\frac{t
  f'\left(\vartheta t\right) Y_t}{f\left(\vartheta
  t\right)}+O\left(\varepsilon \right)\\ 
&=-\frac{tb f'\left(\vartheta
  t\right)}{\varepsilon}\int_{0}^{t}e^{-\int_{s}^{t}q\left(\vartheta
  ,v\right){\rm d}v }\, r_s\,{\rm d}s+\frac{t f'\left(\vartheta t\right)
  Y_t}{f\left(\vartheta t\right)}\\ 
&\qquad \quad -\frac{tb f'\left(\vartheta
  t\right)}{\varepsilon}\int_{0}^{t}e^{-\frac{\left(t-s\right)bf_t}{\varepsilon
} } Y_s\,{\rm d}s+O\left(\varepsilon \right).
\end{align*}
Consider the integral
\begin{align*}
R_\varepsilon \left(t\right)&=\frac{b f\left(\vartheta
  t\right) }{  \varepsilon}\int_{0}^{t}e^{-\frac{\left(t-s\right)bf_t}{\varepsilon
} } \left[Y_s-Y_t\right]\,{\rm d}s\\
&=\int_{0}^{\frac{tbf_t}{\varepsilon }}e^{-u}\left[Y_{t-\frac{u\varepsilon }{bf_t}}-Y_t\right]{\rm d}u=\sqrt{\frac{b\varepsilon}{f_t }}\int_{0}^{\frac{tbf_t}{\varepsilon
}}e^{-u}y_{t,\varepsilon }\left(u\right){\rm d}u,
\end{align*}
where we put $s=t-\frac{u\varepsilon }{bf_t}$ and 
\begin{align*}
y_{t,\varepsilon }\left(u\right)=\sqrt{ {\frac{f_t}{b\varepsilon }
}} \left[{Y_{t-\frac{u\varepsilon
      }{bf_t}}-Y_t}\right]=V_{t,\varepsilon
}\left(u\right)+o\left(\sqrt{\varepsilon }\right) .
\end{align*}

Therefore for any $t\in \left[t_0,T\right]$ as $\varepsilon \rightarrow 0$ we
have the convergence
\begin{align*}
\sqrt{\frac{f_t}{b\varepsilon }}R_\varepsilon \left(t\right)\Longrightarrow
\xi _t\equiv \int_{0}^{\infty } e^{-u}V_t\left(u\right){\rm
  d}u=\int_{0}^{\infty } e^{-u}{\rm d}V_t\left(u\right) 
\end{align*}
and the random variables $\xi _{t_1},\ldots,\xi _{t_k}  $ are independent for
any $0<t_1<\ldots <t_k<T$. Here $V_t\left(\cdot \right)$ is a Wiener process. 

Further
\begin{align*}
&\frac{tb f'\left(\vartheta
  t\right)}{\varepsilon}\int_{0}^{t}e^{-\int_{s}^{t}q\left(\vartheta
  ,v\right){\rm d}v }\, r_s\,{\rm d}s=\frac{tb f'\left(\vartheta
  t\right)}{\varepsilon}\int_{0}^{t}e^{-\frac{\left(t-s\right)bf_t}{\varepsilon
} }\, r_s\,{\rm d}s \left(1+o\left(1\right)\right)\\
&\qquad =\frac{t f'\left(\vartheta
  t\right)}{f\left(\vartheta t\right)}\int_{0}^{\frac{tbf_t}{\varepsilon
  }}e^{-u}\, r_{t-\frac{u\varepsilon }{bf_t}}\,{\rm
    d}u\left(1+o\left(1\right)\right)=\frac{t f'\left(\vartheta 
  t\right)}{f\left(\vartheta t\right)} r_{t
  }\left(1+o\left(1\right)\right).
\end{align*}
Finally we obtain the second presentation \eqref{rep2}: 
\begin{align*}
k_{t,\varepsilon }=-\frac{t f'\left(\vartheta 
  t\right)}{f\left(\vartheta t\right)} r_{t
  }\left(1+o\left(1\right)\right)-\frac{tf'\left(\vartheta
  t\right)}{f\left(\vartheta t\right)} \sqrt{\frac{b\varepsilon
  }{f\left(\vartheta t\right)}} \xi _{t,\varepsilon } \left(1+o\left(1\right)\right)
\end{align*}

From the representations \eqref{rep1} and \eqref{rep2} it follows that 
\begin{align*}
tf'\left(\vartheta t\right)r_t+f\left(\vartheta t\right)k_t=tf'\left(\vartheta
t\right)\sqrt{\frac{b\,\varepsilon 
}{f\left(\vartheta t\right)}}\;\xi
_{t,\varepsilon } \left(1+o\left(1\right)\right)
\end{align*}
and
\begin{align*}
\varepsilon {\rm I}_\varepsilon \left(\vartheta
\right)=b\int_{0}^{T}t^2f'\left(\vartheta t\right)^2\xi
_{t,\varepsilon }^2{\rm d}t +o\left(1\right).
\end{align*}
Of course, $\xi _{t,\varepsilon }, t\in (0,T]$ has no limit process and the
  limit in distribution of each $\xi _{t,\varepsilon }$ is  Gaussian
  random variable  $\xi
  _{t}\sim {\cal N}\left(0,\frac{1}{2}\right) $. The  set  $\xi
  _{t}, t\in (0,T] $ is just a family of independent random variables. 
Let us denote
\begin{align*}
J_\varepsilon\left(\vartheta \right) ={b}\int_{0}^{T}t^2f'\left(\vartheta 
t\right)^2\xi _{t,\varepsilon }^2{\rm d}t.
\end{align*}

 We    have the following properties
\begin{align*}
\lim_{\varepsilon \rightarrow 0}\Ex_{\vartheta } J_\varepsilon
\left(\vartheta  \right)=\frac{{b}}{2}\int_{0}^{T}t^2f'\left(\vartheta 
t\right)^2{\rm d}t\equiv {\rm I}_0 \left(\vartheta  \right),\quad 
\lim_{\varepsilon \rightarrow 0}\Ex_{\vartheta } {J}_\varepsilon
\left(\vartheta  \right)^2= {\rm I}_0 \left(\vartheta  \right)^2,
\end{align*}
which imply that
\begin{align}
\label{qq}
J_\varepsilon
\left(\vartheta  \right)\longrightarrow {\rm I} \left(\vartheta  \right).
\end{align}

{\bf Remark.} Note that the integral 
\begin{align*}
\int_{0}^{T}t^2f'\left(\vartheta t\right)^2\xi
_{t}^2{\rm d}t
\end{align*}
does not exist and the limit \eqref{qq} can be explained as follows. The
Gaussian processes $\xi _{t,\varepsilon }, t\in \left[0,T\right],\varepsilon >0 $ are
continuous and the integral $J_\varepsilon
\left(\vartheta  \right) $   can be well approximated by the sum
\begin{align*}
S_{n,\varepsilon }=\frac{bT}{n}\sum_{j=1}^{n}t_{t_j}^2f'\left(\vartheta
{t_j}\right)^2\xi_{t_j,\varepsilon } ^2,\qquad \quad t_j=\frac{jT}{n}.
\end{align*}
Then we have the first limit ($\varepsilon \rightarrow 0$)
\begin{align*}
S_{n,\varepsilon }\Longrightarrow S_{n
}=\frac{bT}{n}\sum_{j=1}^{n}t_{t_j}^2f'\left(\vartheta {t_j}\right)^2\xi_{t_j}
^2 .
\end{align*}
The second limit ($n\rightarrow \infty $) by the law of large numbers is 
\begin{align*}
S_{n }\longrightarrow  {\rm I}_0 \left(\vartheta  \right)=
\frac{b}{2}\int_{0}^{T}t^2f'\left(\vartheta t\right)^2{\rm d}t. 
\end{align*}
Indeed, we have
\begin{align*}
\Ex_\vartheta S_{n }&=\frac{bT}{2n}\sum_{j=1}^{n}t_{t_j}^2f'\left(\vartheta
   {t_j}\right)^2\longrightarrow \frac{b}{2}\int_{0}^{T}t^2f'\left(\vartheta
   t\right)^2{\rm d}t,\\
 \Ex_\vartheta \left(S_{n    }-\Ex_\vartheta S_{n
 }\right)^2&=\frac{b^2T^2}{n^2}\sum_{j=1}^{n} 
\sum_{i=1}^{n}t_{t_j}^2t_{t_i}^2f'\left(\vartheta
   {t_j}\right)^2f'\left(\vartheta {t_i}\right)^2\Ex_\vartheta\left(\xi_{t_j}
^2-1\right)\left(\xi_{t_i}^2-1\right)\\
&= \frac{Cb^2T^2}{n^2}\sum_{j=1}^{n}
t_{t_j}^4f'\left(\vartheta
   {t_j}\right)^4\leq \frac{ Cb^2T^6K^4}{n}\longrightarrow 0.
\end{align*}

\bigskip

Let us introduce the family of measures $\left\{\Pb_\vartheta^\varepsilon
,\vartheta \in\Theta \right\}$, where $\Pb_\vartheta^\varepsilon $ is the
measure induced in the space of continuous on $\left[0,T\right]$ functions by the
observations $X^T$ satisfying \eqref{a} and define the normalized likelihood ratio
\begin{align*}
Z_\varepsilon \left(u\right)=\frac{V\left(\vartheta
  +{\sqrt{\varepsilon }{u}},X^T\right)}{V\left(\vartheta
  ,X^T\right)} ,\qquad u\in\UU_\varepsilon =\left(\frac{\alpha -\vartheta
    }{\sqrt{\varepsilon }}, \frac{\beta  -\vartheta
    }{\sqrt{\varepsilon }}\right).
\end{align*}

Recall that a statistical experiment is considered as regular in Le Cam's sense if the
corresponding family of measures $\left\{\Pb_\vartheta^\varepsilon
,\vartheta \in\Theta \right\}$ is locally asymptotically normal (LAN)
\cite{IH81}, \cite{H04}. The
studied in the present work model of observations is regular in this  sense. 

\begin{lemma}
\label{L4}
The family of measures $\left\{\Pb_\vartheta^\varepsilon  ,\vartheta
\in\Theta \right\}$ is LAN, i.e., 
we have the representation
\begin{align*}
\ln Z_\varepsilon \left(u\right)=u\Delta _\varepsilon \left(\vartheta
,X^T\right) -\frac{u^2}{2}{\rm I} \left(\vartheta  \right)+\rho _\varepsilon ,
\end{align*}
where $\rho _\varepsilon \rightarrow 0 $,
\begin{align*}
\Delta _\varepsilon
\left(\vartheta,X^T\right)&=\frac{1}{\sqrt{\varepsilon
}}\int_{0}^{T}\left[tf'\left(\vartheta 
  t\right)m\left(\vartheta ,t\right)-f\left(\vartheta
  t\right)\dot m\left(\vartheta,t\right) \right] {\rm d}\bar W_t\\
&=\sqrt{b}\int_{0}^{T}tf'\left(\vartheta
t\right)\xi _{t,\varepsilon }  {\rm d}\bar W_t
\left(1+o\left(1\right)\right)\Longrightarrow  {\cal N}\left(0,{\rm I}_0
\left(\vartheta  \right) \right) .
\end{align*}
\end{lemma}
{\bf Proof.} We have
\begin{align*}
\ln Z_\varepsilon \left(u\right)&=\int_{0}^{T}\frac{M\left(\vartheta
  +\sqrt{\varepsilon }u,t\right)-M\left(\vartheta
   ,t\right)}{{\varepsilon }} \;{\rm d}\bar W_t\\
&\qquad \qquad -\int_{0}^{T}\frac{\left(M\left(\vartheta
  +\sqrt{\varepsilon }u,t\right)-M\left(\vartheta
   ,t\right)\right)^2}{{2\varepsilon^2 }} \; {\rm d}t,
\end{align*}
where
\begin{align*}
M\left(\vartheta  +\sqrt{\varepsilon }u,t\right)-M\left(\vartheta
,t\right)=u\sqrt{\varepsilon } \dot M\left(\vartheta
,t\right)+\frac{u^2{\varepsilon }}{2} \ddot M(\tilde \vartheta
,t). 
\end{align*}
We have the relations
\begin{align*}
&\int_{0}^{T}\frac{\left(M\left(\vartheta
  +\sqrt{\varepsilon }u,t\right)-M\left(\vartheta
   ,t\right)\right)^2}{{\varepsilon^2 }} \; {\rm d}t=\frac{u^2}{\varepsilon }\int_{0}^{T}\dot M\left(\vartheta
   ,t\right)^2 \; {\rm d}t\left(1+o\left(1\right)\right)\\
&\qquad \quad \qquad =u^2\varepsilon  {\rm I}_\varepsilon \left(\vartheta
  \right)\left(1+o\left(1\right)\right) \longrightarrow u^2\;  {\rm
    I}_0\left(\vartheta \right). 
\end{align*}
The asymptotic normality of $\Delta_\varepsilon  \left(\vartheta,X^T \right)$
follows from the central limit theorem for stochastic integrals (see, e.g.,
\cite{Kut94}, Lemma 1.8). 

Let us verify the consistency of the MLE $\hat\vartheta _\varepsilon
$. Consider the log-likelihood ratio
\begin{align*}
\varepsilon \ln \frac{V\left(\vartheta ,X^T\right)}{V\left(\vartheta_0
  ,X^T\right)}& =\int_{0}^{T}{\left[M\left(\vartheta
  ,t\right)-M\left(\vartheta_0
   ,t\right)\right]} \;{\rm d}\bar W_t\\
&\qquad  -\int_{0}^{T}\frac{\left(M\left(\vartheta
  ,t\right)-M\left(\vartheta_0
   ,t\right)\right)^2}{{2\varepsilon }} \; {\rm d}t,
\end{align*}
where we denoted by $\vartheta _0$ the true value. We have to show that the
first integral tends to zero and the second integral tends to a deterministic
function $G\left(\vartheta ,\vartheta _0\right)  $, which has a unique minimum
and the point $\vartheta =\vartheta _0$. 

\begin{lemma}
\label{L5}
We have the convergence
\begin{align}
\label{con}
\varepsilon \ln \frac{V\left(\vartheta ,X^T\right)}{V\left(\vartheta_0
  ,X^T\right)}\longrightarrow - b\int_{0}^{T}\frac{
  \left[f\left(\vartheta t\right)-f\left(\vartheta
    _0t\right)\right]^2}{{4f\left(\vartheta t\right)}} {\rm d}t\equiv -
G\left(\vartheta ,\vartheta _0\right) 
\end{align}
\end{lemma}
{\bf Proof.}
It will be convenient to work with the Kalman filter for the stochastic
process $Z_t=f\left(\vartheta t\right)Y_t$. This leads us to the system of
equations
\begin{align*}
{\rm d}X_t&=Z_t{\rm d}t+\varepsilon {\rm d}W_t,\qquad X_0=0,\\
{\rm d}Z_t&=A\left(\vartheta t\right)Z_t{\rm d}t+bf\left(\vartheta
t\right){\rm d}V_t, \quad Z_0=f\left(0\right)y_0,\qquad 0\leq t\leq T. 
\end{align*}
Here $A\left(\vartheta t\right)=\vartheta \frac{f'\left(\vartheta
  t\right)}{f\left(\vartheta t\right)}-a$. The corresponding filtration
equations are
\begin{align*}
{\rm d}M\left(\vartheta ,t\right)&=A\left(\vartheta t\right)M\left(\vartheta
,t\right){\rm d}t +\frac{\Gamma \left(\vartheta ,t\right)}{\varepsilon
  ^2}\left[{\rm d}X_t- M\left(\vartheta
,t\right){\rm d}t\right],\\
\frac{\partial \Gamma \left(\vartheta ,t\right)}{\partial t}&=2A\left(\vartheta
t\right) \Gamma \left(\vartheta ,t\right)-\frac{\Gamma \left(\vartheta
  ,t\right)^2}{\varepsilon ^2} +b^2f\left(\vartheta t\right)^2,\qquad 0\leq t\leq T,
\end{align*}
with the initial values $M\left(\vartheta ,0\right)=f\left(0\right)y_0 $ and
$\Gamma \left(\vartheta ,0\right)=0$. Here $\Gamma \left(\vartheta
  ,t\right)=\Ex_\vartheta \left(M\left(\vartheta ,t\right)-f\left(\vartheta t\right)Y_t\right)^2 $.  Using the same arguments as above we
obtain a similar to \eqref{asym} approximation 
$$
\Gamma \left(\vartheta
,t\right)=\varepsilon \,b\,f\left(\vartheta t\right)\left(1+O\left(\varepsilon
\right)\right) =\varepsilon \Gamma_* \left(\vartheta
,t\right).
$$ 
If we write the same equations for $M\left(\vartheta _0,t\right)$ and $\Gamma
\left(\vartheta_0 ,t\right) $ and take the difference
$R\left(t\right)=M\left(\vartheta ,t\right)-M\left(\vartheta _0,t\right)  $,
then we obtain the equation for $R\left(t\right)$:
\begin{align*}
{\rm d}R\left(t\right)=B_\varepsilon \left(t\right)R\left(t\right){\rm
  d}t+\left[A\left(\vartheta t\right)-A\left(\vartheta_0 t\right) \right]
M\left(\vartheta _0,t\right){\rm d}t +\delta \left(t\right) {\rm d}\bar W_t, 
\end{align*}
where 
$$ 
B_\varepsilon \left(t\right)=A\left(\vartheta t\right)-\frac{\Gamma_*
  \left(\vartheta ,t\right) }{\varepsilon },\qquad \delta
\left(t\right)=\Gamma_* \left(\vartheta ,t\right)-\Gamma_* \left(\vartheta_0
,t\right).
$$
 The solution of this equation is
\begin{align*}
R\left(t\right)=\int_{0}^{t}e^{\int_{s}^{t}B_\varepsilon \left(v\right){\rm d}v}
\left\{\left[A\left(\vartheta s\right)-A\left(\vartheta_0 s\right) \right] 
M\left(\vartheta _0,s\right){\rm d}s+ \delta \left(s\right) {\rm d}\bar W_s
\right\}. 
\end{align*}
 
For the first integral we have the asymptotics
\begin{align*}
&\int_{0}^{t}e^{\int_{s}^{t}B_\varepsilon \left(v\right){\rm d}v}
\left[A\left(\vartheta s\right)-A\left(\vartheta_0 s\right) \right] 
M\left(\vartheta _0,s\right){\rm d}s\\
&\qquad \qquad \qquad =\varepsilon \frac{\left[A\left(\vartheta
  t\right)-A\left(\vartheta_0 t\right) \right] }{\Gamma_*
  \left(\vartheta ,t\right) } Y_t\left(1+O\left(\sqrt{\varepsilon
}\right)\right). 
\end{align*}
The second integral is of order $\sqrt{\varepsilon }$ because
\begin{align*}
\Ex_{\vartheta _0}\left(\int_{0}^{t}e^{\int_{s}^{t}B_\varepsilon \left(v\right){\rm d}v}
 \delta \left(s\right) {\rm d}\bar W_s
 \right)^2&=\int_{0}^{t}e^{2\int_{s}^{t}B_\varepsilon \left(v\right){\rm d}v} 
 \delta \left(s\right)^2 {\rm d}s\\
&=\frac{\varepsilon \,\delta
   \left(t\right)^2}{2\Gamma_* 
  \left(\vartheta ,t\right) } \left(1+o\left(1\right)\right).
\end{align*}
Further
\begin{align*}
&\int_{0}^{t}e^{\int_{s}^{t}B_\varepsilon \left(v\right){\rm d}v} \delta
\left(s\right) {\rm d}\bar W_s 
 =\frac{\sqrt{\varepsilon }\delta
  \left(t\right)}{\sqrt{\Gamma_* \left(\vartheta ,t\right)} }
\int_{0}^{\frac{t\Gamma_* 
    \left(\vartheta ,t\right) }{\varepsilon }} e^{-u}{\rm d}\bar
W_{t,\varepsilon }\left(u\right)\left(1+o\left(1\right)\right)\\
&\qquad \quad =\frac{\sqrt{\varepsilon }\delta
  \left(t\right)\bar\zeta _{t,\varepsilon }}{\sqrt{\Gamma_* \left(\vartheta
    ,t\right)} }\left(1+o\left(1\right)\right) =\frac{\sqrt{b}
  \left[f\left(\vartheta t\right)-f\left(\vartheta
    _0t\right)\right]}{\sqrt{f\left(\vartheta t\right)}} \bar\zeta
_{t,\varepsilon }\left(1+o\left(1\right)\right) 
\end{align*}
where $\bar W_{t,\varepsilon }\left(u\right) $ is a Wiener process and
$\bar\zeta _{t,\varepsilon } $ is Gaussian random variable. As it was shown
above  the variables $\bar\zeta _{t_1,\varepsilon },\ldots,\bar\zeta
_{t_k,\varepsilon } $ converge in distribution to the independent
i.i.d. random variables $\bar\zeta _{t_1 }, \ldots,\bar\zeta
_{t_k} $, $\bar\zeta _{t }\sim {\cal N}\left(0,1\right) $.

All these allow us  to write
\begin{align*}
\int_{0}^{T}\frac{\left(M\left(\vartheta
  ,t\right)-M\left(\vartheta_0
   ,t\right)\right)^2}{{\varepsilon }} \; {\rm d}t&=b\int_{0}^{T}\frac{
  \left[f\left(\vartheta t\right)-f\left(\vartheta
    _0t\right)\right]^2}{{f\left(\vartheta t\right)}} \bar\zeta
_{t,\varepsilon }^2{\rm d}t\left(1+o\left(1\right)\right)\\
& \longrightarrow  b\int_{0}^{T}\frac{
  \left[f\left(\vartheta t\right)-f\left(\vartheta
    _0t\right)\right]^2}{{2f\left(\vartheta t\right)}} {\rm d}t\equiv
2G\left(\vartheta ,\vartheta _0\right) 
\end{align*}
and
\begin{align*}
\int_{0}^{T}{\left[M\left(\vartheta
  ,t\right)-M\left(\vartheta_0
   ,t\right)\right]} \; {\rm d}\bar W_t\longrightarrow 0.
\end{align*}
It can be shown that all estimates  can be done uniformly in $\vartheta
\in\Theta $. 

Note that the function $G\left(\vartheta ,\vartheta _0\right) $ has a unique
minimum at the point $\vartheta =\vartheta _0$. Moreover, according to
\cite{IH81}, Lemma 3.5.3 we have the estimate
\begin{align*}
G\left(\vartheta ,\vartheta _0\right) \geq c\left|\vartheta -\vartheta _0\right|^2.
\end{align*}
The uniform in $\vartheta $ convergence \eqref{con} provides us the consistency
of the MLE. Recall that the MLE satisfies the equation
\begin{align*}
\dot V\left(\hat\vartheta_\varepsilon  ,X^T\right)&=\int_{0}^{T} 
\frac{\dot M(\hat \vartheta _\varepsilon,t )}{\varepsilon }{\rm d}\bar W_t-\int_{0}^{T} \frac{\dot
M(\hat \vartheta _\varepsilon,t )}{\varepsilon ^2}\left[M(\hat \vartheta
_\varepsilon,t )-M\left(  \vartheta
_0,t \right)\right]{\rm d}t\\
& =\int_{0}^{T} \frac{\dot
M(\hat \vartheta _\varepsilon,t )}{\varepsilon }{\rm d}\bar W_t-  \frac{\left(\hat \vartheta
_\varepsilon-\vartheta _0\right) }{\varepsilon ^2} \int_{0}^{T} \dot
M(\hat \vartheta _\varepsilon,t )  \dot M(\tilde \vartheta
_\varepsilon,t ) {\rm d}t =0.
\end{align*}
Therefore using the consistency of $\hat \vartheta
_\varepsilon$ we can write
\begin{align*}
\frac{\hat \vartheta
_\varepsilon-\vartheta _0}{\sqrt{\varepsilon }}=\frac{\varepsilon ^{-1/2}\int_{0}^{T} \dot
M( \vartheta _0,t )\,{\rm d}\bar W_t }{\varepsilon ^{-1}\int_{0}^{T} \dot
M( \vartheta _0,t )^2 \,  {\rm d}t}\left(1+o\left(1\right)\right).
\end{align*}
Recall that the limit  \eqref{qq} provides us the convergence
\begin{align*}
\varepsilon ^{-1/2}\int_{0}^{T} \dot
M( \vartheta _0,t )\,{\rm d}\bar W_t\Longrightarrow {\cal N}\left(0, {\rm
  I}_0\left(\vartheta _0\right)\right), \quad \varepsilon ^{-1}\int_{0}^{T} \dot
M( \vartheta _0,t )^2 \,  {\rm d}t\longrightarrow {\rm
  I}_0\left(\vartheta _0\right).
\end{align*}
Hence the asymptotic normality \eqref{an} is proved.

\section{Discussions}

The results on frequency estimation by the observations of stationary Gaussian
process $Y_t$ satisfying linear equation in the presence of WGN are valid for  much
more general models of inhomogeneous processes $Y_t$ and   (smooth) functions
$f\left(\vartheta ,t\right)$. We did not use the periodicity of
$f\left(\vartheta ,t\right)=f\left(\vartheta t\right)$ and the relation 
\begin{align*}
\lim_{\varepsilon \rightarrow 0} \int_{0}^{T}\dot M\left(\vartheta ,t\right)^2{\rm d}t=0
\end{align*}
holds for inhomogeneous processes $Y_t$ too. Recall that we took this model of
observations just for the comparison of the properties of estimators for
different limits. 

 Recall that if the smooth signal $f\left(\vartheta t\right)$ is observed in
 the WGN (say, $Y_t\equiv 1$ in \eqref{a}) and $T\rightarrow \infty $, then
 the rate of convergence of the MLE $\hat\vartheta _T$ is $T^{3/2}$
 \cite{IH1}:
\begin{align*}
T^{3/2}\left(\hat\vartheta _T-\vartheta \right)\Longrightarrow {\cal
  N}\left(0,{\rm I}_*\left(\vartheta \right)^{-1}\right)
\end{align*}
with some ${\rm I}_*\left(\vartheta \right)>0 $.
In the case $\sigma =b =\varepsilon\rightarrow 0$ studied in \cite{ChK1} we have
\begin{align*}
\frac{\hat\vartheta _\varepsilon -\vartheta }{\varepsilon }\Longrightarrow {\cal
  N}\left(0,{\rm I}\left(\vartheta \right)^{-1}\right).
\end{align*} 
and the limit variance ${\rm I}\left(\vartheta \right)^{-1}$ for large $T$ is
of order $T^{-3}$, i.e.; if we consider the {\it second limit } $T\rightarrow
\infty $, then the normalization formally can be written as follows
\begin{align*}
\frac{T^{3/2}\left(\hat\vartheta _\varepsilon -\vartheta\right) }{\varepsilon
}\Longrightarrow {\cal N}.
\end{align*} 
In our case $\sigma =\varepsilon\rightarrow 0 $  and $b>0$ fixed we have
\begin{align*}
\frac{\hat\vartheta _\varepsilon -\vartheta }{\sqrt{\varepsilon \,b}
}\Longrightarrow {\cal 
  N}\left(0,\tilde {\rm I}\left(\vartheta \right)^{-1}\right),
\end{align*}
where 
\begin{align*}
\tilde {\rm I}\left(\vartheta \right)=\int_{0}^{T}t^2f'\left(\vartheta
t\right)^2{\rm d}t .
\end{align*}
Therefore if $b\rightarrow 0$ (second limit), then we can write formaly the
normalization $\sqrt{\varepsilon b}$. If we consider now the {\it third limit}
$T\rightarrow \infty $, then we can wait that
\begin{align*}
  \frac{T^{3/2}\left(\hat\vartheta _\varepsilon -\vartheta\right)
  }{\sqrt{\varepsilon \,b} }\Longrightarrow {\cal 
  N}.
\end{align*}

 We have to note that the calculation of the MLE for the model of partially
 observed linear system is of extreme computational complexity because to
 calculate it we have to know the solutions of the filtration equations
 \eqref{kf}, \eqref{eR} for all $\vartheta \in \Theta $. To realize an
 estimator asymptotically equivalent to the MLE and much more easy calculated
 we can use the Multi-step MLE approach developed in \cite{Kut17},
 \cite{KhK17}. In the case of periodic signal in WGN the similar One-step MLE
 approach was realized in \cite{G84}. The most interesting case of frequency
 estimation for the model of observation \eqref{a},\eqref{b} is $T\rightarrow
 \infty $ and it will be considered in our next work. For periodic diffusion
 processes the similar  problems of frequency estimation were considered in
 \cite{HK11}, \cite{HK12}

{\bf Acknowledgment.} This work was done under partial financial support of
the grant of RSF number 14-49-00079.

\end{document}